\documentclass{article}

\usepackage[latin1]{inputenc}
\usepackage{amsmath, amsthm, amsfonts}
\usepackage[english]{babel}
\usepackage{graphics}
\usepackage{graphicx}
\usepackage{amsfonts}
\usepackage{amssymb}
\usepackage{subfigure}
\usepackage{ragged2e}

\theoremstyle{definition}

\theoremstyle{remark}

\newtheorem{defi}{Definition}[section]

\newtheorem*{Proof}{Proof}

\newtheorem*{teo2}{Keywords}

\newtheorem{lema}[defi]{Lemma}
\newtheorem{teo}[defi]{Theorem}

\theoremstyle{definition}
\newtheorem{ejem}[defi]{Example}
\theoremstyle{remark}
\newtheorem*{remark}{Remark}

\title{Ribaucour-type surfaces }
\author{Milton Javier Cardenas Mendez \\ Instituto de Matemática e Estatística \\
Universidade Federal de Goiás, 74001-970 Goiânia, GO, Brazil\\ miltonjcardenas@ufg.br \and  Armardo Mauro Vasquez Corro\\ Instituto de Matemática e Estatística \\ Universidade Federal de Goiás, 74001-970 Goiânia, GO, Brazil\\corro@ufg.br }
\date{\today}

\begin{document}
\maketitle

\begin{abstract}
In this work we define the Ribaucour-type surfaces (in short, RT-surfaces). These surfaces satisfy a relationship similar to the Ribaucour surfaces that are related to the Élie Cartan problem. This class furnishes what seems to be the first examples of pairs of noncongruent
surfaces in Euclidean space such that, under a diffeomorphism, lines of curvatures
are preserved and principal curvatures are switched. We show that every compact and connected RT-surface is a sphere with center at the origin.
We obtain present a Weierstrass type representation for RT-surfaces with prescribed Gauss map which depends on two holomorphic functions. We give explicit examples of RT-surfaces. Also, we use this representation to classify the RT-surfaces of rotation.
\end{abstract}

\begin{teo2}
Ribaucour surfaces, generalized Weingarten surfaces, prescribed normal Gauss map, Weierstrass type representation.
\end{teo2}

\textbf{Introduction}: An oriented surface $\Sigma \subseteq \mathbb{R}^3$ is called a Weingarten surface if it exists
a differentiable relation $W$ between the mean curvatures $H$ and Gaussian curvatures $K$ of $\Sigma$, such that
$W(H,K) = 0$, when  function $W$ is linear, that is, $a+bH +cK = 0$ for $ a, b$ and $ c$
constants, the surfaces are called linear Weingarten surfaces.

Examples of linear Weingarten surfaces are constant Gaussian curvature surfaces ($c \neq 0$ and $b = 0$) and  constant mean curvature surfaces ($b \neq 0$ and $c = 0$). Several authors have studied these classes of surfaces, see \cite{a}, \cite{b} and \cite{c}, among others.

Let $\Sigma\subset \mathbb{R}^3$ be an oriented surface with normal Gauss map $N$, functions $\Psi,\Lambda :\Sigma\rightarrow \mathbb{R}^3$ given by $ \Psi(p) = \langle p , N(p) \rangle$ and $ \Lambda(p) = \langle p , p \rangle $, $ p \in \Sigma$, where $\langle,\rangle$ denotes the Euclidean scalar product in $\mathbb{R}^3$, are called support function
and quadratic distance function, respectively. Geometrically, $\Psi(p)$ measures the signed distance from the origin to $T_pM$ and $\Lambda(p)$ measures the squared distance from the origin to $p$. Let $p\in \Sigma $, a sphere with center $p+\frac{H(p)}{K(p)} N(p)$ and radius $\frac{H(p)}{K(p) }$ is called the middle sphere.

In 1888, Appell \cite{f} studied a class of surfaces oriented in $ \mathbb{R}^3 $ associated with area-preserving sphere transformations. Later, Ferreira and Roitman \cite{d} showed that these surfaces satisfy the Weingarten relation, $H + \Psi K = 0 $.

In 1907, Tzitzeica \cite{o} studied hyperbolic surfaces oriented so that there is a nonzero constant $c \in \mathbb{R} $ for which $ K + c^2 \Psi^4 = 0 $.

In \cite{a}, the authors motivated by the works \cite{f} and \cite{o} defined  generalized Weingarten surfaces as surfaces that satisfy a relation of the form $A + BH + CK = 0$, where $A, B, C : \Sigma \rightarrow \mathbb{R}$ are differentiable functions that do not depend on the parameterization of $\Sigma$.
  In particular, they studied the class of surfaces that satisfy the relation $2\Psi H + \Lambda K = 0$. Called Special Generalized Weingartem Surfaces depending on the support function and the distance function (in short, EDSGW-surfaces), these surfaces have the geometric property that all medium spheres pass through the origin. The authors obtained a Weierstrass-like representation of EDSGW-surfaces depending on two homomorphic functions. In \cite{b}, the authors classified isothermic EDSGW-surfaces in relation to the third fundamental shape parameterized by plane curvature lines. Also in \cite{z}, it is shown that EDSGW-surfaces are in correspondence with the surface class in $\mathbb{S}^2 \times \mathbb{R}$ where the Gaussian curvature $ K$ and the extrinsic curvature $K_E$ satisfy $K=2K_E$.

Martínez and Roitman, in \cite{c} showed what appears to be the first example found for the second case of the problem posed by Élie Cartan in his classic book about external differential systems and their applications to Differential Geometry. Such examples are given by a class of Weingarten surfaces that satisfy the relation $2 \Psi H + (1+ \Lambda) K = 0 $, Ribaucour surface calls, these surfaces have the geometric property that all the medial spheres intercept a fixed sphere along a large circle.

In \cite{pp}, the authors define a surface class called  Ribaucour surface of harmonic type (in short, HR-surface) if it satisfies
$2\Psi H + (c e^{2\mu} +\Lambda)K = 0$, where $c$ is a nonzero real constant, $\mu$ a harmonic function with respect to the third
fundamental form. These surfaces generalize  Ribaucour surfaces studied in \cite{c}.

Motivated by \cite{a},\cite{b}, \cite{c} and \cite{d}, we define  Ribaucour-type surfaces (RT-surfaces) which have the geometric property for every $ p \in \Sigma $ a sphere of center $p + (\frac{H (p)} {K (p) } + \frac {\Psi (p)} {2}) N (p) $ e radius $ \frac {H (p)} {K (p)} + \frac {\Psi (p)} {2} $ passes through the origin, in this case $2 \Psi H + (\Lambda + \Psi^2) K = 0 $ is satisfied.
We obtain a Weierstrass-type representation for  RT-surfaces with prescribed normal Gaussian application depending on two holomorphic functions. Using this representation we classify rotating RT-surfaces. Furthermore, we show that a compact and connected RT-surface is a sphere.

\section{Preliminaries}
In this section we fix the notation used in this work, $ \Sigma $ a surface on $ \mathbb{R}^3 $, $ N $ its normal Gaussian map, and $ U $ an open subset of $ \mathbb{R}^2 $.

Let $ X: U \subset \mathbb{R}^2 \rightarrow \Sigma$, a parameterization of a surface $ \Sigma $ and $ N: U \subset \mathbb{R}^2 \rightarrow \mathbb{R}^{3} $, normal Gaussian map. Considering $ \{X_{,1}, X_{,2}, N \} $ as a base of $ \mathbb{R}^{3} $,  where $ X_{,i} (q) = \frac{\partial X} {\partial u_i} (q) $, $ 1 \leq i, j \leq 2 $ further we can write  vector $ X_{, ij} $, $1 \leq i, j \leq 2 $, as
\begin{equation*}
  X_{,ij}= \sum_ {k = 1}^ {2} \widetilde{\Gamma}_ {ij}^k X_{,k} + b_{ij} N
\end{equation*}
The $ \widetilde{\Gamma}_ {ij}^k $ coefficients are called Christoffel symbols.

\begin{defi}
\label{def1}
Let $ X $ be a local parameterization of $ \Sigma $ with map of Gauss $ N $, matrix $ W = (W_{ij})$, such that
\begin{equation*}
N_{,i} = \sum_{j = 1}^{2} W_{ij} X_j, \hspace{0.3cm} 1 \leq i \leq 2
\end{equation*}
is called the Weingarten matrix of $\Sigma$.
\end{defi}

\begin{lema}
\label{le1}
Let $N$ be the normal Gaussian map given by (\ref{N1})  such that the metric $ L_{ij} =\langle N_{, i}, N_{, j} \rangle $ is Euclidean conformal.  Christoffel's symbols for metric $ L_ {ij}$ are given by
\begin{equation*}
\Gamma_{ij}^k=0, \hspace{0.3cm} \Gamma_{ii}^i=\frac{L_{ii,i}}{2L_{ii}},\hspace{0.3cm} \Gamma_{ij}^i=\frac{L_{ii,j}}{2L_{ii}},\hspace{0.3cm} \Gamma_{ii}^i=\frac{-L_{ii,j}}{2L_{jj}}
\end{equation*}
For $i,j$ e $k$  different.
\end{lema}

The next result was obtained by Roitman and Ferreira  \cite{d}.
\begin{teo}
\label{te1}
Let $ \Sigma \subset \mathbb{R}^{n + 1} $ be an orientable hypersurface and $ N: \Sigma \rightarrow \mathbb{S}^n $ normal Gauss application with non-zero Gauss-kroncker curvature in every point.  Let $ U \subset \Sigma $ be a neighborhood of $ p_0 $ such that $ N: U \rightarrow N(U) = V \subset \mathbb{S}^n $ invertible and $ h(q) = \langle q, N^{- 1}(q) \rangle, q \in V$, then
\begin{equation*}
X(q)=\nabla_Lh(q)+h(q)N(q)
\end{equation*}
\end{teo}
\begin{remark}
Using the above equation we have functions
\begin{equation}
\label{qu}
\Lambda(q)=\langle X(q), X(q)\rangle=|\nabla_Lh(q)|^2+h(q)^2
\end{equation}
\begin{equation}
\label{su}
\Psi(q)=\langle X(q), N(q)\rangle=h(q)
\end{equation}
called quadratic distance function  and support function, respectively.
\end{remark}

\begin{remark}
\label{holo}
 Let the inner product be defined by $ \langle, \rangle:\mathbb{C} \times \mathbb{C} \rightarrow \mathbb{R}, \langle f,g \rangle= f_1g_1+f_2g_2$, where $f= f_1+if_2$ and $g=g_1+ig_2$ are holomorphic functions. If $f,g,h:U\subset\mathbb{C}\rightarrow\mathbb{C}$ are holomorphic functions, then
\begin{equation}\label{holo}
     \begin{split}
         &\langle f,g\rangle_{,1}=\langle f',g\rangle+\langle f,g'\rangle  \\
                    &\langle f,g\rangle_{,2}=\langle if',g\rangle+\langle f,ig'\rangle\\
             &\langle fh,g\rangle=\langle f,\overline{h} g\rangle\\
    &f =\langle 1,g\rangle+i\langle i,f\rangle\\
 & g_{,1}=g', g_{,2}=ig'\\
  \end{split}
\end{equation}
where  $\langle f,g\rangle_{,1}=\frac{\partial\langle f(z),g(z)\rangle}{\partial u_1} $
 \end{remark}

\section{RT-surfaces}
Motivated by the works \cite{a},\cite{b}, \cite{c} and \cite{d}, we will begin the study of Ribaucour-type surfaces  and will call them RT-surfaces. In addition to presenting some examples, will provide a Weierstrass representation depending on two holomorphic functions for surfaces in this class and characterize the case where such surfaces are of rotation.

\begin{teo}
\label{te2}
Let $\Sigma \subset \mathbb{R}^3 $, an orientable surface with non-zero Gauss-Kronecker  curvature. Then there is  a differentiable  function  $h:U \rightarrow \mathbb{R}$ and  $g$ a holomorphic function  such that normal Gauss map  is given by
 \begin{equation}
\label{N1}
  N=\frac{( 2g(u),1-|g(u)|^2)}{1+|g(u)|^2}
\end{equation}
the coefficients of the $III$ fundamental form are
\begin{equation}
\label{simb}
L_{ij}=\frac{4 |g'|^ 2 \delta _{ij} }{(1+|g|^ 2)^ 2}
\end{equation}
 $\Sigma $ is locally parameterized by
\begin{equation}
\label{1}
X(u)=\sum_{j=1}^{2} \frac{h(u)_{,j}}{L_{jj}}N(u)_{,j}+h(u)N(u)
\end{equation}
In this case   $h(u)=\langle X(u),N(u)\rangle$ is the support function. Furthermore, the Weingarten matrix is given by
$ W = V^{-1} $ where
\begin{equation}
\label{2}
V_{ij}=\frac{1}{L_{ij}}\left( h_{,ij}-\sum_{k}^{n}h_{,k}\Gamma_{ij}^{k}+hL_{ij}\delta_{ij}\right)
\end{equation}
where $\Gamma_{ij}^{k}$ are Christoffel's symbols  of $N$ and the fundamental forms $ I $ and $ II $ of $ X $, in local coordinates, are given by

 \begin{equation*}
 I=\langle X_{,i},X_{,j}\rangle=\sum_{k=1}^{n} V_{ik}V_{jk}L_{kk},\quad II=\langle X_{,i},N_{,j}\rangle=V_{ij}L_{jj}.
 \end{equation*}
\end{teo}

\begin{Proof}
In theorem (\ref{te1}) we can choose a local parameterization of the sphere given by (\ref{N1}) with a metric given by (\ref{simb}) and $\Sigma$ is locally parameterized by (\ref{1}), calculating their derivatives we get
\begin{equation*}
X_{,i}=\left(\frac{h_{,ii}}{L_{ii}}-\sum_{j=1}^{n}\frac{h_{,j}}{L_{ii}} \Gamma_{ii}^j +h \right)N_{,i}+\sum_{\substack {j=1\\ j\neq i}}^{n}\left(\frac{h_{,ij}}{L_{jj}}-\frac{h_{,j}}{L_{jj}}\Gamma_{ij}^j-\frac{h_{,i}}{L_{jj}}\Gamma_{ij}^i \right)N_{,j}
\end{equation*}
Considering matrix $V=(V_{ij}),\hspace{0.3cm} 1\leq i,j\leq n$, given by $(\ref{2})$, therefore
\begin{equation}
\label{14}
X_{,i}=\sum_{j=1}^{n}V_{ij}N_{,j}
\end{equation}
In search of the Weingarten matrix we have to $N_{,i}=\sum_{j=1}^{n}V^{-1}_{ij}X_{,j}$, by definition \ref{def1} we have to $W=V^{-1}$. To obtain the coefficients of the fundamental forms, we use (\ref{14}), therefore
\begin{equation*}
I =\langle X_{,i},X_{,j}\rangle=\langle\sum_{k=1}^{n}V_{ik}N_{,k},\sum_{m=1}^{n}V_{jm}N_{,m}\rangle
=\sum_{k,m=1}^{n}V_{ik} V_{jm}\langle N_{,k}N_{,m}\rangle=\sum_{k=1}^{n}V_{ik} V_{jk}L_{kk}
\end{equation*}
\begin{equation*}
II =\langle X_{,i},N_{,j}\rangle=\langle\sum_{k=1}^{n}V_{ik}N_{,k},N_{,j}\rangle=\sum_{k=1}^{n}V_{ik}\langle N_{,k},N_{,j}\rangle=V_{ij}L_{jj}
\end{equation*}
\end{Proof}

\begin{defi}
A $\Sigma $ surface is called Ribaucour-type surfaces (RT-surface) such that for each $ p \in \Sigma $ the center sphere $ p + (\frac{H(p)}{K (p )} + \frac{\Psi(p)}{2}) N (p) $ and radius $ \frac{H(p)}{K(p)} + \frac{\Psi(p)}{2 } $ go through the origin, in this case $ \Sigma $ satisfies the following generalized Weingarten relation
\begin{equation*}
  2\Psi H+(\Lambda+\Psi^2)K=0
\end{equation*}
for all $p\in \Sigma$.
\end{defi}

\begin{lema}
\label{ss}
Let $ \Sigma $ be a Riemann surface and $ X: \Sigma \rightarrow \mathbb{R}^3 $ an immersion such that the Gauss-kronecker curvature is non-zero, under the conditions of the theorem \ref{te2}  then $ X(\Sigma) $ is a SS-surface if and only if
\begin{equation}\label{lap}
 h \triangle h- | \nabla h |^ 2 = 0
\end{equation}
\end{lema}
\begin{proof}
By theorem (\ref{te2}), we have $W=V^{-1}$, let $\sigma_i$ be the eigenvalues  of $V$ and  $\lambda_i$ the eigenvalues  of $W$, then  $\sigma_i=\frac{1}{\lambda_i}$. Using  this fact  and the  expressions in lemma \ref{le1}, (\ref{simb}), (\ref{1}) and (\ref{2}), we get
\begin{equation*}
   \frac{\triangle h}{L_{11}}+2h=V_{11}+V_{22}=\frac{-2H}{K}
\end{equation*}
using equations (\ref{qu}) and (\ref{su}),  then
\begin{equation*}
 \left(\frac{h\triangle h-|\nabla h|^2}{L_{11}}\right)K+ 2H\Psi+ (\Lambda+\Psi^2)K=0,\quad \forall p\in \Sigma
 \end{equation*}
and the result follows.
\end{proof}

For RT-surfaces with Gaussian curvature $ K \neq 0 $ we will present a complete characterization through pairs of holomorphic functions. This representation will allow to classify all RT-surfaces of rotation. Before that, we will need the following lemma.

\begin{lema}
\label{lemav}
Consider holomorphic functions $f,g: U\subset \mathbb{C} \rightarrow \mathbb{C} $, with $ g'\neq 0 $. Taking  local parameters $ z = u_1 + iu_2 \in \mathbb{C} $ and $ h = e^{\langle 1, f \rangle} $. So elements of the matrix $V$ given by (\ref{2}) in terms of $f$ and $g$ are given by
\begin{equation}
\label{vij}
\begin{split}
   & V_{11} =\frac{T^2h}{4|g'|^2}[\langle 1,f'\rangle^2-\langle 1,\xi \rangle] + h \\
     & V_{12}=V_{21}=\frac{T^2h}{4|g'|^2}\langle i,\xi -\frac{f'^2}{2} \rangle  \\
     &   V_{22}=\frac{T^2h}{4|g'|^2}[\langle 1,if'\rangle^2+\langle 1,\xi \rangle] + h
\end{split}
\end{equation}
where $T=1+ | g | ^ 2 $ and  $\xi =f' \left(\frac{g''}{g'}-\frac{2}{T}g'\overline{g}\right)-f''$.  Furthermore,
\begin{equation}\label{33}
detV=\frac{h^2T^4}{16|g'|^4}\langle \xi,f'^2-\xi\rangle+ \frac{h^2T^2 |f'|^2}{4|g'|^2}+h^2
\end{equation}
\end{lema}

\begin{Proof}
With $h=e^{\langle 1,f\rangle}$, deriving and using (\ref{holo}) we get
\begin{equation}
\label{h}
\begin{split}
&h_{,1}=e^{\langle 1,f\rangle}\langle 1,f'\rangle, \quad h_{,11} = e^{\langle 1,f\rangle}(\langle 1,f'\rangle^2+\langle 1,f''\rangle),\quad h_{,2}=e^{\langle 1,f\rangle}\langle 1,if'\rangle,\\
     &h_{,22} = e^{\langle 1,f\rangle}(\langle 1,if'\rangle^2-\langle 1,f''\rangle),\qquad h_{,12} = e^{\langle 1,f\rangle}(\langle 1,if'\rangle\langle 1,f'\rangle+\langle 1,if''\rangle).
\end{split}
\end{equation}
Since $T=1+|g|^2$, by lemma \ref{le1} and (\ref{simb}), we obtain
\begin{equation*}
   \Gamma_{11}^1 =\frac{T\langle g',g''\rangle-2|g'|^2\langle g,g'\rangle }{T|g'|^2}=\Gamma_{12}^2, \quad \Gamma_{22}^2=\frac{T\langle g',ig''\rangle-2|g'|^2\langle g,ig'\rangle }{T|g'|^2}=\Gamma_{21}^1
\end{equation*}
\begin{equation*}
\Gamma_{11}^2=\frac{2|g'|^2\langle g,ig'\rangle -T\langle g',ig''\rangle}{T|g'|^2},\quad\Gamma_{22}^1=\frac{2|g'|^2\langle g,g'\rangle-T\langle g',g''\rangle }{T|g'|^2}.
\end{equation*}
Using the above expressions, (\ref{2}) and (\ref{h}), we obtain (\ref{vij}), so we have
\begin{equation*}
\begin{split}
   detV &=\left(\frac{hT^2}{4|g'|^2}[\langle 1,f'\rangle^2-\langle 1,\xi \rangle] + h\right) \left( \frac{hT^2}{4|g'|^2}[\langle 1,if'\rangle^2+\langle 1,\xi \rangle] + h\right)\\
     &-\left(\frac{hT^2}{4|g'|^2}[\langle 1,if'\rangle\langle 1,f'\rangle+\langle i,\xi \rangle]\right)^2
\end{split}
\end{equation*}
therefore we obtain(\ref{33}).
\end{Proof}

\begin{teo}
\label{te3}
Let $ \Sigma $ be a Riemann surface and $ X: \Sigma \rightarrow \mathbb{R}^3 $ an immersion such that the Gauss-kronecker curvature is non-zero. \\
Then $ X (\Sigma) $ is a SS-surface if and only if there are holomorphic functions $ f, g: U\subset \mathbb{C} \rightarrow \mathbb{C} $, where $U$ is a simply connected open and  $ g' \neq 0 $, such that $ X (\Sigma) $ is locally parameterized by
\begin{equation}
\label{p1}
X=\frac{e^{\langle 1,f\rangle}}{2|g'|^2}(Tg'\bar{f'}-2g\langle g',gf'\rangle , -2\langle g',gf'\rangle) + e^{\langle 1,f\rangle}\frac{(2g,2-T)}{T}
\end{equation}
With normal application of Gauss N given by (\ref{N1}), the regularity condition is given by
\begin{equation*}
T^4 \langle \xi,f'^2-\xi \rangle  +4T^2|f'|^2|g'|^2+ 16|g'|^4 \neq 0
\end{equation*}
The coefficients of the first and second fundamental forms of $ X $ have the following expressions
\begin{equation*}
\begin{split}
E&= h^2\left(\frac{T^2}{4|g'|^2}\left(A_1^2+\langle i,\xi-\frac{f'^2}{2} \rangle^2\right)+2A_1+\frac{4|g'|^2}{T^2}\right)\\
F&=h^2\left(\frac{T^2|f'|^2}{4|g'|^2} +2\right)\langle i,\xi-\frac{f'^2}{2} \rangle\\
   G&=h^2\left(\frac{T^2}{4|g'|^2}\left(A_2^2+\langle i,\xi-\frac{f'^2}{2} \rangle^2\right)+2A_2+\frac{4|g'|^2}{T^2}\right)\\
e&=hA_1 + \frac{4h|g'|^2}{T^2},\quad f=h\langle i,\xi-\frac{f'^2}{2} \rangle, \quad g=hA_2 +\frac{4h|g'|^2}{T^2}\end{split}
\end{equation*}
Where $\xi =f' \left(\frac{g''}{g'}-\frac{2}{T}g'\overline{g}\right)-f''$, $A_1=\langle 1,f'\rangle^2-\langle 1,\xi \rangle$ e $A_2=\langle 1,if'\rangle^2-\langle 1,\xi \rangle$
\end{teo}

\begin{Proof}
Using (\ref{lap} ), we can assume  without loss of generality  $ h = e^{\mu} $ where $ \mu: \Sigma \rightarrow \mathbb{R} $ is a differentiable function, in this case,
\begin{equation*}
h \triangle h- | \nabla h |^ 2 = e ^ {2\mu} \triangle\mu
\end{equation*}
Now $h$ satisfies (\ref{lap}), if and only if $h=e^{\langle 1,f\rangle}$, where $f$ is a holomorphic function. As $N$ is given by (\ref{N1}), deriving we have
\begin{equation*}
   N_{,1}=\frac{2}{T^2}(T g'-2g\langle g',g\rangle , -2\langle g',g \rangle)
\end{equation*}
  \begin{equation*}
   N_{,2}=\frac{2}{T^2}(T ig'-2g\langle g,ig'\rangle , -2\langle g,ig' \rangle  )
\end{equation*}
using these expressions and (\ref{simb}) in (\ref{1}),we have (\ref{p1}).
Therefore, using (\ref{holo}),(\ref{14}) and (\ref{vij}) we obtain the coefficients of the first and second fundamental forms given in the statement of the theorem.
\end{Proof}
\begin{ejem}
A sphere of center at the origin  and radius $r>0 $   is RT-surface, in fact, using theorem \ref{te3} and taking $h(u)=r$, by equations (\ref{1}) and (\ref{lap}) the result follows.
\end{ejem}

\begin{teo}
Let $ \Sigma $ be a compact, connected SS-surface, then $ \Sigma $ is a sphere with center at the origin.
\end{teo}
\begin{Proof}
  Let $ \Sigma $ be a compact  then there is a $ E $ sphere with center at the origin of radius $r>0$, such that $ \Sigma $  is contained in the closed ball with center at the origin and radius $r$, $B[0,r]$ and a point $ p \in E \cap \Sigma $, such that $ T_pE = T_p\Sigma $.
Let $ h: \Sigma \rightarrow \mathbb{R} $, where  $ h $ function support of $ \Sigma $. We know that $ h(q) \leq |q| $, for every $ q \in \Sigma $, then
\begin{equation*}
h(q)\leq |q| \leq r  \quad \text{for all} \quad q\in \Sigma
\end{equation*}
$\Sigma$ is a RT-surface and by theorem \ref{te3} and (\ref{1}), there is a  parameterized locally around $p$ such that
support function $h:U\subset \mathbb{R}^2 \rightarrow \Sigma $,
given by $h(u)=e^{\mu(u)}$ $u \in U$ , with $\mu$ harmonic then
\begin{equation*}
e^{\mu}\leq r \Rightarrow \mu\leq \ln(r)
\end{equation*}
Since $\mu$ is harmonic, by the maximum principle, $\mu(u)= \ln(r)$ in $U$, later $h(u)=r$ in $U$, for an argument of compactness and connectedness of $\Sigma$, we conclude   $h(q)=r$ for all $q\in \Sigma $, therefore $\Sigma$  is a sphere  with center at the origin.
\end{Proof}

For some holomorphic functions $f$ and $g$ we show some examples of (\ref{p1}).

\begin{figure}[htbp]
\begin{center}
\subfigure{\includegraphics[width=32mm]{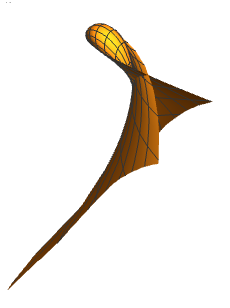}}
\subfigure{\includegraphics[width=32mm]{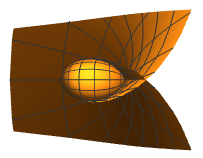}}
\subfigure{\includegraphics[width=32mm]{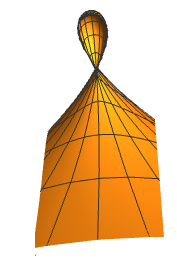}}
\end{center}
\caption{$f(z)=g(z)=z=u_1+iu_2$}
\label{lash3}
\end{figure}

\begin{figure}[htbp]
\begin{center}
\subfigure{\includegraphics[width=32mm]{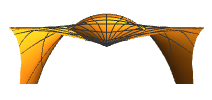}}
\subfigure{\includegraphics[width=32mm]{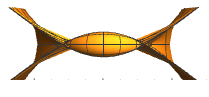}}
\subfigure{\includegraphics[width=32mm]{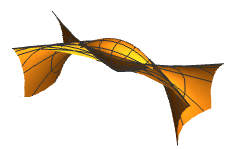}}
\end{center}
\caption{$f(z)=z^2,\hspace{0.2cm} g(z)=z$}
\label{lash3}
\end{figure}



\begin{figure}[htbp]
\begin{center}
\subfigure{\includegraphics[width=32mm]{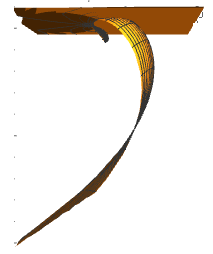}}
\subfigure{\includegraphics[width=32mm]{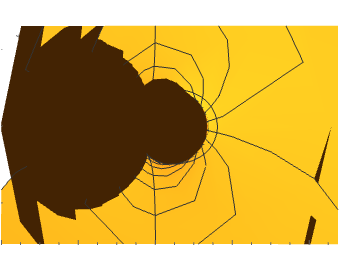}}
\subfigure{\includegraphics[width=32mm]{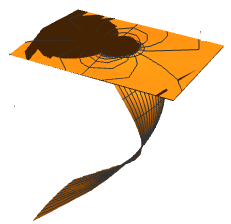}}
\end{center}
\caption{$f(z)=z,\hspace{0.2cm} g(z)=z^2$}
\label{lash3}
\end{figure}

\newpage
The following theorem characterizes the rotating RT-surfaces.
\begin{teo}
Let $ \Sigma $ be a connected RT-surface. Since $ \Sigma $ is of rotation if and only if there are constants $ a, b \in \mathbb{R} $, such that $ \Sigma $ can be locally parameterized by
\begin{equation}\label{3.5}
X_{a,b}(u_1,u_2)=\frac{e^{au_1+b}}{1+e^{2u_1}} \left(M(u_1)\cos(u_2),M(u_1)\sin(u_2),N(u_1)\right)
\end{equation}
where
\begin{equation}\label{3.51}
M(u_1)=  \frac{ a(e^{-u_1}-e^{3u_1})+ 4 e^{u_1}}{2},\quad N(u_1)=1-e^{2u_1}-a(1+e^{2u_1})
\end{equation}
\end{teo}

\begin{Proof}
Let $\Sigma$ be locally parameterized by  (\ref{p1}) with normal application of Gauss  $N$ given by (\ref{N1}) and where $ f, g $ are holomorphic functions. As $ \Sigma $ is of rotation if and only if $ g (w) = w $ and $ h (w) = J (| w |^ 2) $, $ w \in \mathbb{C} $, for some $ J $ differentiable function. Changing parameters $w=e^z, z=u_1+iu_2\in \mathbb{C}$, we have $ g (z) = e^z $ and $ h (z) = J(e^{2u_1}) $.
Consequently, $ h,_2 = 0 $, remembering that $ h = e^{\langle 1, f \rangle} $ and $ h,_2 = e^{\langle 1, f \rangle} \langle 1, if'\rangle = 0 $, then $ \langle 1, if' \rangle = 0 $, so
\begin{equation*}
f(z)=az+z_0,\hspace{.2cm}g(z)=e^z, \hspace{.2cm} h(z)=e^{au_1+b},\hspace{.2cm}z=u_1+iu_2,z_0=b+ic\in \mathbb{C}
\end{equation*}
using these expressions in (\ref{p1}), we have
\begin{equation*}
   X=e^{au_1+b}\left(\left(\frac{a(e^{-u_1}-e^{u_1})}{2}+\frac{2e^{u_1}}{1+e^{2u_1}}\right)(\cos(u_2)+i\sin(u_2)),\frac{1-e^{2u_1}}{1+e^{2u_1}}-a\right)
\end{equation*}
Therefore  $\Sigma$ can be parameterized locally by (\ref{3.5}) and (\ref{3.51}).
\end{Proof}

Using theorem \ref{te3}, we have that
\begin{equation*}
EG-F^2=\frac{(ae^{2u_1}(e^{2u_1}-4)-a)^2 (ae^{4u_1}(a+1)+2e^{2u_1}(a^2+2)+a(a-1))^2}{16 (e^{2u_1}+1)^4 e^{-4u_1(a-1)-4b}}
\end{equation*}
Thereby, $X_{ab}$ is regular if and only if $a=0$. If $a\neq 0$, the expression from above  vanishes for

\begin{equation*}
u_1=\left\{ \begin{array}{lcc}
             \frac{1}{2}\ln( \frac{2+ \sqrt{4+ a^2}}{a}) &   se  & a>0 \\
             \\ \frac{1}{2}\ln( \frac{2-\sqrt{4+ a^2}}{a}) &   se  & a<0 \\
             \\ \frac{1}{2}\ln( \frac{-(2+a^2)+ \sqrt{4+ 5a^2}}{a^2+a}) &   se  & 0<a<1\\
             \\ \frac{1}{2}\ln( \frac{-(2+a^2)- \sqrt{4+ 5a^2}}{a^2+a}) &   se  & -1<a<0
             \end{array}
   \right.
\end{equation*}
In figures (4), (5) and (6) we present examples of RT-surfaces of rotation. We consider only the case where $b = 0$.

\begin{figure}[h]
\begin{minipage}[b]{0.3\linewidth}
\centering
\includegraphics[width=\linewidth]{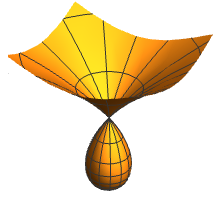}
\caption{$a=-1$}
\label{fig:figura1}
\end{minipage}
\hspace{0.3cm}
\begin{minipage}[b]{0.22\linewidth}
\centering
\includegraphics[width=\linewidth]{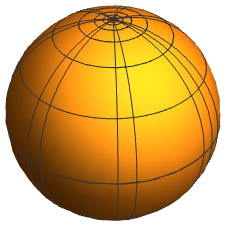}
\caption{$a=0$}
\label{fig:figura1}
\end{minipage}
\hspace{0.3cm}
\begin{minipage}[b]{0.49\linewidth}
\centering
\includegraphics[width=\linewidth]{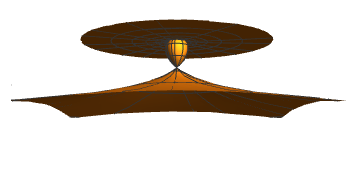}
\caption{$a=\frac{1}{2}$}
\label{fig:figura2}
\end{minipage}
\end{figure}

\renewcommand{\refname}{Bibliografía}








\end{document}